\documentclass{interact}
\usepackage{bm}
\def\<{\langle}
\def\>{\rangle}
\usepackage{upgreek}

\usepackage{amsmath, amsthm, amscd, amsfonts, amssymb, graphicx, color}
\usepackage[bookmarksnumbered, colorlinks, plainpages]{hyperref}
\hypersetup{colorlinks=true,linkcolor=red, anchorcolor=green, citecolor=red, urlcolor=red, filecolor=magenta, pdftoolbar=true}

\newtheorem{theorem}{Theorem}[section]
\newtheorem{lemma}[theorem]{Lemma}
\newtheorem{proposition}[theorem]{Proposition}
\newtheorem{corollary}[theorem]{Corollary}
\theoremstyle{definition}

\newtheorem{example}[theorem]{Example}
\newtheorem{remark}[theorem]{Remark}

\theoremstyle{definition}

\newcommand{\Hc}{\mathbb{H}_c}
\newcommand{\C}{\mathbb{C}}
\newcommand{\R}{\mathbb{R}}
\newcommand{\Fix}{\mathrm{Fix}}
\newcommand{\Nm}{\mathrm{N}}
\newcommand{\tr}{\operatorname{tr}}
\newcommand{\bq}[1]{\overline{#1}}
\def\dss{\displaystyle}
\theoremstyle{remark}

\numberwithin{equation}{section}

\begin{document}
	\title{ Binormal Toeplitz Matrices over Commutative Quaternions}
	
	\author{
		\name{Muhammad Ahsan Khan}
		\affil{Department of Mathematics, University of Kotli Azad Jammu $\&$ Kashmir, Kotli 11100, Azad Jammu $\&$ Kashmir, Pakistan, Email ahsankhan388@hotmail.com}
	}
	
	\maketitle

	\begin{abstract}
		In this paper, Toeplitz matrices over the algebra of commutative quaternions are studied through binormality. Since this algebra is a product of two copies of the complex numbers, the first conjugation reduces to the complex case, while the other two are equivalent and give a coupled condition on a pair of complex Toeplitz matrices. We describe binormal Toeplitz matrices completely in sizes two and three and for bidiagonal matrices of every size. We also show that the shift property of Campbell holds for the first conjugation when the scalar is not a zero divisor, but fails for the other two.
	\end{abstract}
	
	
	\begin{amscode}
		13E10,	47A05, 15B05 
	\end{amscode}	
	
	\section{Introduction}
	Hamilton discovered the quaternions in the year 1843, and for for a long time the word quaternion meant a number system that gives up commutativity of multiplication. Almost fifty years later Segre observed that one can do the opposite: keep commutativity and give up the property that every non-zero number can be inverted \cite{Se}. The resulting four-dimensional algebra is called the algebra of commutative quaternions; it is also known as the algebra of bicomplex numbers or of commutative reduced biquaternions \cite{Pr}. It contains the complex numbers, it has three imaginary units, and it contains non-zero numbers whose product is zero, the so-called zero divisors.
	
	There are several reasons for working with commutative quaternions rather than with the Hamilton quaternionss. First, commutativity makes matrix computations behave as over a commutative ring, with determinants and characteristic polynomials, and the eigenvalues of a square matrix are ordinary complex numbers \cite{KT1,KAT}; for Hamilton quaternions left and right eigenvalues must be distinguished and no determinant with the usual properties exists. Second, it has three different conjugations, hence three different adjoints of a matrix. Commutative quaternions have been proposed for colour image and signal processing \cite{Pei}, and Toeplitz matrices are the classical models of linear filters and of covariance matrices of stationary sequences \cite{GS,BS,Gr}. Our aim here is structural: to understand what changes when the scalars are commutative quaternions.

	The property we studied is binormality. A matrixx is normal when it commutes with its adjoint, and binormal when the two products of the matrix with its adjoint, taken in the two orders, commute. The class was introduced by Campbell \cite{C1,C2}; it contains all normal matrices and many non-normal ones, and for complex Toeplitz matrices of small size binormal matrices were characterized by Ko, Kwon and Lee \cite{KKL}. Toeplitz matrices with quaternion entries, in particular with commutative quaternion entries, have been studied only little \cite{Kl,KT1,KAT,KY2,KK}, and binormality has not been considered in this setting.

	Because commutative quaternions are a product of two copies of the complex numbers, a matrix over them is a pair of complex matrices. For the first conjugation, binormality of the matrix is binormality of both members, so the complex theory applies and only the bookkeeping caused by the zero divisors is new. For the second and third conjugations the adjoint exchanges the two members, and binormality becomes a coupled condition on a pair of complex matrices. These two conjugations are equivalent equivalent  through the automorphism that exchanges two of the imaginary units, so there are only two genuinely different cases. For them the natural quadratic form is indefinite, so binormality is understood algebraically and not through absolute values as in Campbell's work; the meaning of this notion is discussed in Section 2.

	The main results of the paper are the following. First, we give a complete and explicit description of binormal Toeplitz matrices of size three over the commutative quaternions: such a matrix is binormal exactly when two explicit vectors, whose entries are quadratic in the entries of the matrix, are parallel (Theorem~\ref{thm:n3full}). To our knowledge this is the first explicit description of binormal Toeplitz matrices beyond size two. It is important because normal Toeplitz matrices fall into two simple types, whereas binormal ones are governed by a proportionality condition; the description is uniform in all three conjugations through the reduction to pairs of complex Toeplitz matrices. Second, the proof exhibits the structure behind the equations: the commutator of the two factors is skew-persymmetric and cyclic, which turns binormality into a statement about a combination of the identity and the square of that commutator. Third, we show that Campbell's shift property, which says that a matrix is normal as soon as it and its translate by a non-zero scalar are binormal, fails for the second and third conjugations in the strongest possible way: there is a matrix whose translate by every scalar is binormal while the matrix is not normal (Example~\ref{ex:genfail}); for the first conjugation the property holds exactly when the scalar is not a zero divisor (Theorem~\ref{thm:shift1} and Example~\ref{ex:shiftfail}).

	The paper is organized as follows. Section 2 recalls the algebra, reduces the third conjugation to the second, and contains the pair formulation, a discussion of binormality for an indefinite form, and the general commutator identity. Section 3 contains the Toeplitz results: sizes two and three, bidiagonal and monomial matrices, and the shift property.

	\section{Commutative Quaternions, Matrices and Binormality}
	
	The algebra of commutative quaternions is
	\[
	\Hc=\{q=t+\mathbf{i}x+\mathbf{j}y+\mathbf{k}z:\ t,x,y,z\in\R\},
	\]
	with the multiplication determined by $\R$-bilinearity and
	\[
	\mathbf{i}^2=\mathbf{k}^2=-1,\qquad \mathbf{j}^2=1,\qquad \mathbf{ij}=\mathbf{ji}=\mathbf{k},\qquad \mathbf{jk}=\mathbf{kj}=\mathbf{i},\qquad \mathbf{ki}=\mathbf{ik}=-\mathbf{j} .
	\]
	All products of two basis elements are symmetric in the two factors, so the multiplication is commutative; associativity will follow from Lemma~\ref{lem:split}, which identifies $\Hc$ with a product of fields. We denote by $\C=\{t+\mathbf{i}x\}$ the copy of the complex numbers inside $\Hc$. Since $\mathbf{k}=\mathbf{ij}=\mathbf{ji}$, every $q\in\Hc$ can be written uniquely as
	\[
	q=c_1+\mathbf{j}c_2,\qquad c_1=t+\mathbf{i}x\in\C,\quad c_2=y+\mathbf{i}z\in\C ,
	\]
	because $\mathbf{j}c_2=\mathbf{j}y+\mathbf{j}\mathbf{i}z=\mathbf{j}y+\mathbf{k}z$.
	
	\begin{lemma}\label{lem:split}
		Put $e_+=\tfrac12(1+\mathbf{j})$ and $e_-=\tfrac12(1-\mathbf{j})$.$\beta\beta$
		\begin{enumerate}
			\item[(a)] $e_+^2=e_+$, $e_-^2=e_-$, $e_+e_-=0$, $e_++e_-=1$ and $\mathbf{j}=e_+-e_-$.
			\item[(b)] The map $\psi:\Hc\to\C\times\C$, $\psi(q)=(a,b)$, where
			\[
			a=c_1+c_2=(t+y)+\mathbf{i}(x+z),\qquad b=c_1-c_2=(t-y)+\mathbf{i}(x-z),
			\]
			is an isomorphism of $\R$-algebras, whose inverse is $(a,b)\mapsto ae_++be_-$. In particular $\Hc$ is commutative and associative.
		\end{enumerate}
	\end{lemma}
	
	\begin{proof}
		(a) is a direct computation using $\mathbf{j}^2=1$. (b) For $q=c_1+\mathbf{j}c_2$ we have $q=c_1(e_++e_-)+c_2(e_+-e_-)=(c_1+c_2)e_++(c_1-c_2)e_-$, and conversely $ae_++be_-=\frac{a+b}2+\mathbf{j}\frac{a-b}2$; so $\psi$ is a bijection with the stated inverse, and it is $\R$-linear. It is multiplicative because $e_\pm^2=e_\pm$, $e_+e_-=0$ and $\C$ commutes with $\mathbf{j}$: $(ae_++be_-)(a'e_++b'e_-)=aa'e_++bb'e_-$. A bijective multiplicative linear map onto the commutative associative algebra $\C\times\C$ transfers these properties to $\Hc$.
	\end{proof}
	
	We call $a$ and $b$ the \emph{coordinates} of $q$ and write $q=(a,b)$. For any $q=t+\mathbf{i}x+\mathbf{j}y+\mathbf{k}z,,,,$
	\[
	q^{(1)}=t-\mathbf{i}x+\mathbf{j}y-\mathbf{k}z,\qquad q^{(2)}=t+\mathbf{i}x-\mathbf{j}y-\mathbf{k}z,\qquad q^{(3)}=t-\mathbf{i}x-\mathbf{j}y+\mathbf{k}z
	\]
	are the three non-trivial involutive automorphisms of $\Hc$. 
	
	\begin{lemma}\label{lem:conj}
		In the coordinates of Lemma~\ref{lem:split},
		\[
		q^{(1)}=(\bar a,\bar b),\qquad q^{(2)}=(b,a),\qquad q^{(3)}=(\bar b,\bar a).
		\]
		Consequently each $q\mapsto q^{(\kappa)}$ is a ring automorphism with $(q^{(\kappa)})^{(\kappa)}=q$, and $(q^{(1)})^{(2)}=q^{(3)}$ and cyclically. The fixed algebras are
		\[
		\Fix_1=\{t+\mathbf{j}y\}\cong\R\times\R,\qquad \Fix_2=\{t+\mathbf{i}x\}=\C,\qquad \Fix_3=\{t+\mathbf{k}z\}\cong\C .
		\]
		Thus $\Fix_1$ has zero divisors, whereas $\Fix_2$ and $\Fix_3$ are fields.
	\end{lemma}
	
	\begin{proof}
		Write $q=c_1+\mathbf{j}c_2$. Reversing the signs of $x,z$ gives $q^{(1)}=\bar c_1+\mathbf{j}\bar c_2$, with coordinates $(\bar a,\bar b)$; reversing the signs of $y,z$ gives $q^{(2)}=c_1-\mathbf{j}c_2$, with coordinates $(c_1-c_2,c_1+c_2)=(b,a)$; reversing the signs of $x,y$ gives $q^{(3)}=\bar c_1-\mathbf{j}\bar c_2$, with coordinates $(\bar b,\bar a)$. Conjugation of both coordinates, the swap, and their composition are ring automorphisms of $\C\times\C$ of order two, and applying the first and then the second gives $(a,b)\mapsto(\bar b,\bar a)$.

		Fixed points: $q^{(1)}=q\iff a,b\in\R\iff q\in\mathrm{span}_\R\{1,\mathbf{j}\}$ (as $e_\pm\in\mathrm{span}\{1,\mathbf{j}\}$); $q^{(2)}=q\iff a=b\iff q\in\C$; $q^{(3)}=q\iff b=\bar a$, and then with $a=\alpha+\mathbf{i}\beta$ one gets $q=\alpha(e_++e_-)+\mathbf{i}\beta(e_+-e_-)=\alpha+\beta\mathbf{k}$, and $\alpha+\beta\mathbf{k}\mapsto\alpha+\beta\mathbf{i}$ is an isomorphism onto $\C$ because $\mathbf{k}^2=-1$.
	\end{proof}
	
\begin{lemma}\label{lem:tau}
		Let $\tau:\Hc\to\Hc$ be the $\R$-linear map $\tau(t+\mathbf{i}x+\mathbf{j}y+\mathbf{k}z)=t+\mathbf{i}z+\mathbf{j}y+\mathbf{k}x$. Then $\tau$ is an automorphism of $\R$-algebras, in coordinates $\tau(a,b)=(a,\bar b)$, and
		\[
		\tau\big(q^{(2)}\big)=\big(\tau q\big)^{(3)},\qquad \tau\big(q^{(3)}\big)=\big(\tau q\big)^{(2)} .
		\]
		Consequently, applying $\tau$ entrywise, a matrix $T$ is binormal (resp.\ normal) for $\kappa=2$ if and only if $\tau(T)$ is binormal (resp.\ normal) for $\kappa=3$, and $\tau$ preserves the Toeplitz structure.
	\end{lemma}
	
	\begin{proof}
		For $q=c_1+\mathbf{j}c_2$ one has $\tau q=(t+\mathbf{i}z)+\mathbf{j}(y+\mathbf{i}x)$, with coordinates $(a,\bar b)$. The map $(a,b)\mapsto(a,\bar b)$ is an $\R$-algebra automorphism of $\C\times\C$. By Lemma~\ref{lem:conj}, $\tau(q^{(2)})=\tau(b,a)=(b,\bar a)=(a,\bar b)^{(3)}$ and $\tau(q^{(3)})=\tau(\bar b,\bar a)=(\bar b,a)=(a,\bar b)^{(2)}$. Entrywise, $\tau(T^{\dagger_2})=\tau(T)^{\dagger_3}$, and $\tau$ is a ring automorphism of $M_n(\Hc)$, so it preserves products and commutators.
	\end{proof}
	
	By Lemma~\ref{lem:tau} the second and third conjugations are one case: every statement for $\kappa=3$ is the statement for $\kappa=2$ applied to $\tau(T)$. We keep the uniform proofs for $\kappa=2,3$ below, and in the examples the computation for $\kappa=3$ is the computation for $\kappa=2$ of the matrix $\tau(T)$.
	
		From now on a conjugation index $\kappa\in\{1,2,3\}$ is fixed, we write $\bq q=q^{(\kappa)}$ and
	\[
	\Nm(q)=\bq q\,q\in\Fix_\kappa .
	\]
	In coordinates, by Lemma~\ref{lem:conj},
	\begin{equation}\label{eq:N}
		\Nm(q)=(|a|^2,|b|^2)\ \ (\kappa=1),\qquad \Nm(q)=(ab,ab)\ \ (\kappa=2),\qquad \Nm(q)=(\bar ba,\ \bar ab)\ \ (\kappa=3).
	\end{equation}
	
	\begin{lemma}\label{lem:units}
		An element $q=(a,b)$ is invertible $\iff$ $a\ne0$ and $b\ne0$, $\iff$ $q$ is not a zero divisor. The non-zero zero divisors are exactly the non-zero elements of $e_+\Hc$ and of $e_-\Hc$, that is, the elements with exactly one vanishing coordinate. Moreover, $q\,q^{(1)}q^{(2)}q^{(3)}=|a|^2|b|^2$, and for the three imaginary units
		\begin{equation}\label{eq:Nijk}
			\Nm(\mathbf{i})=(1,-1,1),\qquad \Nm(\mathbf{j})=(1,-1,-1),\qquad \Nm(\mathbf{k})=(1,1,-1)\qquad\text{for }\kappa=1,2,3 .
		\end{equation}
	\end{lemma}
	
	\begin{proof}
		In $\C\times\C$ the element $(a,b)$ has inverse $(a^{-1},b^{-1})$ when $ab\ne0$, while $(0,b)(1,0)=0$ shows that $(0,b)$ is a zero divisor, and similarly if $b=0$. By Lemma~\ref{lem:conj}, $q\,q^{(1)}q^{(2)}q^{(3)}=(a,b)(\bar a,\bar b)(b,a)(\bar b,\bar a)=(|a|^2|b|^2,|a|^2|b|^2)$. For \eqref{eq:Nijk} use $\bar{\mathbf{j}}=\mathbf{j},-\mathbf{j},-\mathbf{j}$, $\bar{\mathbf{i}}=-\mathbf{i},\mathbf{i},-\mathbf{i}$, $\bar{\mathbf{k}}=-\mathbf{k},-\mathbf{k},\mathbf{k}$ for $\kappa=1,2,3$ together with $\mathbf{i}^2=\mathbf{k}^2=-1$, $\mathbf{j}^2=1$.
	\end{proof}
	
	Two features of \eqref{eq:N} will be used constantly. First, $\Nm(q)$ lies in $\{t+\mathbf{j}y\}$ only for $\kappa=1$; for $\kappa=2,3$ it lies in the other fixed algebras. Second, for $\kappa=2,3$ the quantity $\Nm(q)$ is not positive in any sense: it can be negative, and it can vanish for $q\ne0$.

	Throughout we designate by $M_n(\Hc)$ the space of all $n\times n$ matrices with entries all from $\Hc$. For $T=(t_{ij})\in M_n(\Hc)$ let $\bq T=(\bq t_{ij})$ and $T^\dagger=(\bq T)^t$ (the adjoint with respect to the chosen conjugation). The matrix $T$ is \emph{normal} if $TT^\dagger=T^\dagger T$, \emph{selfadjoint} if $T=T^\dagger$, and \emph{binormal} if $T^\dagger T$ and $TT^\dagger$ commute. Because entries commute and the conjugation is a ring automorphism, $(XY)^\dagger=Y^\dagger X^\dagger$ and $(X^\dagger)^\dagger=X$, and the diagonal entries of a selfadjoint matrix lie in $\Fix_\kappa$. We write $[X,Y]=XY-YX$.
	
	Applying $\psi$ entrywise, each $T\in M_n(\Hc)$ can be written uniquely as $T=e_+A+e_-B$ with $A,B\in M_n(\C)$; we call $A$ and $B$ the \emph{components} of $T$. Because $e_\pm$ are central orthogonal idempotents, the map $T\mapsto(A,B)$ is a ring isomorphism from $M_n(\Hc)$ onto $M_n(\C)\times M_n(\C)$: indeed $(e_+A+e_-B)(e_+A'+e_-B')=e_+AA'+e_-BB'$, because the $(i,j)$ entry of the left side is the sum over $k$ of products of entries, to which the identity $(ae_++be_-)(a'e_++b'e_-)=aa'e_++bb'e_-$ applies termwise. 
	
	\begin{proposition}\label{prop:red}
		Let $T=e_+A+e_-B$ with $A,B\in M_n(\C)$.
		\begin{enumerate}
			\item[(i)] $\kappa=1$: $T^\dagger=e_+A^*+e_-B^*$. $T$ is normal $\iff$ $A$ and $B$ are normal, and binormal $\iff$ $A$ and $B$ are binormal.
			\item[(ii)] $\kappa=2$: $T^\dagger=e_+B^t+e_-A^t$. $T$ is normal $\iff$ $B^tA=AB^t$, and binormal $\iff$ $[B^tA,AB^t]=0$.
			\item[(iii)] $\kappa=3$: $T^\dagger=e_+B^*+e_-A^*$. $T$ is normal $\iff$ $B^*A=AB^*$, and binormal $\iff$ $[B^*A,AB^*]=0$.
		\end{enumerate}
	\end{proposition}
	
	\begin{proof}
		By Lemma~\ref{lem:conj}, the $(i,j)$ entry $\bq t_{ji}$ of $T^\dagger$ has coordinates $(\bar a_{ji},\bar b_{ji})$, $(b_{ji},a_{ji})$, $(\bar b_{ji},\bar a_{ji})$ for $\kappa=1,2,3$; these are the $(i,j)$ entries of $(A^*,B^*)$, $(B^t,A^t)$, $(B^*,A^*)$. (i) By the ring isomorphism above, $T^\dagger T=e_+A^*A+e_-B^*B$ and $TT^\dagger=e_+AA^*+e_-BB^*$, and $e_+X+e_-Y$, $e_+X'+e_-Y'$ are equal, resp.\ commute, if and only if $X=X'$ and $Y=Y'$, resp.\ $[X,X']=0$ and $[Y,Y']=0$. (ii) Here $T^\dagger T=e_+B^tA+e_-A^tB$ and $TT^\dagger=e_+AB^t+e_-BA^t$; since $A^tB=(B^tA)^t$, $BA^t=(AB^t)^t$ and $[X^t,Y^t]=-[X,Y]^t$, the conditions on the second components are the transposes of those on the first. (iii) The same, with $(B^*A)^*=A^*B$, $(AB^*)^*=BA^*$ and $[X^*,Y^*]=-[X,Y]^*$.
	\end{proof}
	
	Thus for $\kappa=1$ the problem consists of two independent complex problems, while for $\kappa=2,3$ it is a coupled problem for the pair $(A,B^t)$, respectively $(A,B^*)$. Notice that if $A$ and $B$ are Toeplitz, then so are $B^t$ and $B^*$.
	
	\begin{proposition}\label{prop:pair}
		Let $T=e_+A+e_-B$. Put $(A,C)=(A,A^*)$ and $(B,B^*)$ for $\kappa=1$ (two pairs), $(A,C)=(A,B^t)$ for $\kappa=2$, and $(A,C)=(A,B^*)$ for $\kappa=3$. Then $T$ is binormal if and only if $[CA,AC]=0$ for each pair, and normal if and only if $CA=AC$ for each pair. If $T$ is Toeplitz, then $A$ and $C$ are complex Toeplitz matrices; for $\kappa=2,3$ the matrix $C$ is an arbitrary Toeplitz matrix, and $\kappa=1$ is the special case $C=A^*$.
	\end{proposition}
	
	Campbell's definition uses the Hilbert space adjoint, for which $T^*T$ and $TT^*$ are positive and binormality says that $|T|$ and $|T^*|$ commute. For $\kappa=2,3$ the form $\Nm(q)=\bq q$ is indefinite: it can be negative and it can vanish for $q\ne0$. Binormality is therefore defined here algebraically, as the commutation $[T^\dagger T,TT^\dagger]=0$ in the ring $M_n(\Hc)$ with its involution; it agrees with Campbell's notion for $\kappa=1$. By Proposition~\ref{prop:pair} it is exactly the condition $[CA,AC]=0$, with normality $CA=AC$ as the extreme case, which makes it the natural notion. Some phenomena below already occur for complex matrices when the transpose replaces the adjoint, for instance $T^tT=0$ with $T=\begin{pmatrix}1&0\\ \mathbf{i}&0\end{pmatrix}\ne0$, or a negative diagonal entry of $T^tT$; they reflect the indefinite form and not the zero divisors, and we indicate where this applies (Examples~\ref{ex:neg} and~\ref{ex:isotropic}).

	The complex adjoint $\Phi:M_n(\Hc)\to M_{2n}(\C)$ of Kosal and Tosun \cite{KT1} replaces each entry $c_1+\mathbf{j}c_2$ by the block $c_1I+c_2\sigma$, $\sigma=\begin{pmatrix}0&1\\1&0\end{pmatrix}$. Since $\sigma$ has eigenvalues $\pm1$ with fixed orthonormal eigenvectors, $\Phi(T)$ is unitarily similar to $A\oplus B$ for every $\kappa$, and $\Phi(T^{\dagger_1})=\Phi(T)^*$. Hence $\Phi(T)$ is binormal if and only if $A$ and $B$ are, which for $\kappa=1$ is binormality of $T$. For $\kappa=2,3$, however, binormality of $T$ is the coupled condition of Proposition~\ref{prop:red}(ii),(iii) and is not equivalent to binormality of $\Phi(T)$ (see Example~\ref{ex:Phi}).

\begin{lemma}\label{lem:crit}
		Let $T\in M_n(\Hc)$, $U=T^\dagger T$, $V=TT^\dagger$, $\Delta=U-V$ and $M=\Delta U$. Then $[U,V]=M-M^\dagger$. In particular $T$ is binormal $\iff$ $M$ is selfadjoint, i.e.\ $M_{ij}=\overline{M_{ji}}$ for all $i,j$.
	\end{lemma}
	
	\begin{proof}
		Since $V=U-\Delta$, $[U,V]=U(U-\Delta)-(U-\Delta)U=\Delta U-U\Delta$. Both $U$ and $\Delta$ are selfadjoint, so $M^\dagger=U^\dagger\Delta^\dagger=U\Delta$, and $[U,V]=M-M^\dagger$.
	\end{proof}
	
	\section{ Binormal Toeplitz Matrices with Commuting Quaternion Entries}
	
	A Toeplitz matrix over $\Hc$ is a matrix $T=(t_{i-j})_{i,j=1}^n$ with $t_m\in\Hc$ for $-(n-1)\le m\le n-1$, constant along each diagonal. We put $t_m=0$ for $|m|\ge n$. Since the condition $t_{ij}=t_{i-j}$ is $\R$-linear, $T=e_+A+e_-B$ is Toeplitz $\iff$ $A$ and $B$ are Toeplitz. We write $J$ for the exchange matrix with entries $J_{ij}=\delta_{i,n+1-j}$; it is real, symmetric and satisfies $J^2=I$. A Toeplitz matrix is persymmetric, $JT^tJ=T$, because $(JT^tJ)_{ij}=(T^t)_{n+1-i,\,n+1-j}=T_{n+1-j,\,n+1-i}=t_{(n+1-j)-(n+1-i)}=t_{i-j}$. Consequently $T^t=JTJ$, and $T$ is binormal $\iff$ $T^t$ is binormal, since conjugation by the real orthogonal matrix $J$ preserves commutation.
	
	Let $T=(t_{i-j})_{i,j=1}^n$ be a Toeplitz matrix over $\Hc$ and put
	\[
	U=T^\dagger T,\qquad V=TT^\dagger,\qquad \Delta=U-V .
	\]
	Then $T$ is normal $\iff$ $\Delta=0$ and binormal $\iff$ $[U,V]=0$.
	
	\begin{lemma}\label{lem:uv}
		The entries of $U$ and $V$ are
		\[
		u_{ij}=\sum_{k=1}^n\bq t_{k-i}\,t_{k-j},\qquad v_{ij}=\sum_{k=1}^n t_{i-k}\,\bq t_{j-k},
		\]
		and both matrices are selfadjoint: $u_{ji}=\bq u_{ij}$ and $v_{ji}=\bq v_{ij}$. In particular all diagonal entries of $U$, $V$ and $\Delta$ lie in $\Fix_\kappa$.
	\end{lemma}
	
	\begin{proof}
		The formulas follow from $(T^\dagger)_{ik}=\bq t_{k-i}$ and the definition of the product. Conjugating $u_{ij}$ gives $\dss\sum_kt_{k-i}\bq t_{k-j}=u_{ji}$ by commutativity, and similarly for $v$; a diagonal entry equals its own conjugate, i.e.\ lies in $\Fix_\kappa$.
	\end{proof}
	
	\begin{lemma}\label{lem:V}
		For all $i,j$,
		\[
		v_{ij}=u_{n+1-j,\,n+1-i}=\overline{u_{n+1-i,\,n+1-j}},\qquad\text{equivalently}\qquad V=J\,\bq U\,J=J\,U^t\,J .
		\]
	\end{lemma}
	
	\begin{proof}
		In the sum defining $v_{ij}$ substituting $\ell=n+1-k$ and put $a=n+1-i$, $b=n+1-j$: $v_{ij}=\dss\sum_{\ell}t_{\ell-a}\bq t_{\ell-b}=\sum_\ell\bq t_{\ell-b}t_{\ell-a}=u_{ba}=\bq u_{ab}$. Since $U$ is selfadjoint, $\bq U=U^t$, which gives the matrix form.
	\end{proof}
	
	\begin{lemma}\label{lem:Delta}
		$\Delta$ is selfadjoint and skew-persymmetric: $\Delta_{n+1-j,\,n+1-i}=-\Delta_{ij}$, i.e.\ $J\Delta^tJ=-\Delta$. Moreover $\Delta_{ij}=u_{ij}-\overline{u_{n+1-i,\,n+1-j}}$.
	\end{lemma}
	
	\begin{proof}
		By Lemma~\ref{lem:V}, $V=JU^tJ$, hence $U=JV^tJ$ and $J\Delta^tJ=JU^tJ-JV^tJ=V-U=-\Delta$. Selfadjointness follows from that of $U$ and $V$, and the last formula from Lemma~\ref{lem:V}.
	\end{proof}

	\subsection{The $2\times 2$ Case}
	
	Let $T=\begin{pmatrix}t_0&t_{-1}\\t_1&t_0\end{pmatrix}$. By Lemma~\ref{lem:uv},
	\[
	u_{11}=\Nm(t_0)+\Nm(t_1),\qquad u_{22}=\Nm(t_0)+\Nm(t_{-1}),\qquad u_{12}=\bq t_0t_{-1}+\bq t_1t_0,\qquad u_{21}=\bq u_{12},
	\]
	and by Lemma~\ref{lem:V}
	\[
	V=\begin{pmatrix}u_{22}&u_{12}\\ \bq u_{12}&u_{11}\end{pmatrix},\qquad \Delta=\begin{pmatrix}d&0\\0&-d\end{pmatrix},\qquad d:=u_{11}-u_{22}=\Nm(t_1)-\Nm(t_{-1})\in\Fix_\kappa .
	\]
	
	\begin{lemma}\label{lem:claim}
		If $t_0$ is invertible and $u_{12}=0$, then $d=0$.
	\end{lemma}
	
	\begin{proof}
		The equation $u_{12}=0$ says $\bq t_0\,t_{-1}=-\bq t_1\,t_0$. Since $t_0$ is invertible, so is $\bq t_0$, and $t_{-1}=-(\bq t_0)^{-1}\bq t_1t_0$. Conjugating, $\bq t_{-1}=-t_0^{-1}t_1\bq t_0$. Multiplying the two identities and using commutativity, $\Nm(t_{-1})=\bq t_{-1}t_{-1}=t_1\bq t_1=\Nm(t_1)$, i.e.\ $d=0$.
	\end{proof}
	
	\begin{theorem}\label{thm:n2}
		Let $T$ be a $2\times 2$ Toeplitz matrix over $\Hc$ and let $\kappa\in\{1,2,3\}$.
		\begin{enumerate}
			\item[(1)] $T$ is binormal $\iff$ $d\,u_{12}=0$, and normal $\iff$ $d=0$. In particular $d=0$ implies that $T$ is binormal, with no hypothesis on $t_0$; and $t_0=0$ implies $u_{12}=0$, hence that $T$ is binormal.
			\item[(2)] If $t_0$ is not a zero divisor, then $T$ is binormal $\iff$ $T$ is normal.
			\item[(3)] If $\kappa\in\{2,3\}$, then $T$ is binormal $\iff$ $T$ is normal or $t_0=0$.
			\item[(4)] If $\kappa=1$, then $T$ is binormal $\iff$ for each sign $\varepsilon\in\{+,-\}$ the component $e_\varepsilon T$ is normal or has zero diagonal.
		\end{enumerate}
	\end{theorem}
	
	\begin{proof}
		(1) By Lemma~\ref{lem:crit}, $T$ is binormal if and only if $M=\Delta U=\begin{pmatrix}du_{11}&du_{12}\\-d\,\bq u_{12}&-du_{22}\end{pmatrix}$ is selfadjoint. The diagonal entries lie in $\Fix_\kappa$, so the only condition is $M_{12}=\overline{M_{21}}$, i.e.\ $du_{12}=-du_{12}$ (because $\bq d=d$), i.e.\ $du_{12}=0$. Normality means $d=0$. If $d=0$ or $t_0=0$ (then $u_{12}=\bq t_0t_{-1}+\bq t_1t_0=0$), the matrix is binormal.

		(2) Let $t_0$ be invertible and $du_{12}=0$. If $\kappa\in\{2,3\}$ and $d\ne0$, then $d$ is invertible in the field $\Fix_\kappa$, so $u_{12}=0$ and Lemma~\ref{lem:claim} gives $d=0$. If $\kappa=1$, multiply by $e_\varepsilon$: the ring $e_\varepsilon\Hc\cong\C$ is a field stable under conjugation, so $e_\varepsilon d=0$ or $e_\varepsilon u_{12}=0$, and in the second case Lemma~\ref{lem:claim}, applied in $e_\varepsilon\Hc$ (its proof uses only commutativity and the invertibility of $t_0$), gives $e_\varepsilon d=0$. Hence $d=0$. The converse is clear, since normal matrices are binormal.

		(3) Let $\kappa\in\{2,3\}$, let $T$ be binormal and $d\ne0$; we show $t_0=0$. As $\Fix_\kappa$ is a field, $u_{12}=0$. Write $t_m=(a_m,b_m)$. By Lemma~\ref{lem:conj}, $u_{12}=\big(b_0a_{-1}+b_1a_0,\ a_0b_{-1}+a_1b_0\big)$ for $\kappa=2$ and $u_{12}=\big(\bar b_0a_{-1}+\bar b_1a_0,\ \bar a_0b_{-1}+\bar a_1b_0\big)$ for $\kappa=3$, while $d=(\delta,\delta)$ with $\delta=a_1b_1-a_{-1}b_{-1}$, respectively $d=(\delta,\bar\delta)$ with $\delta=\bar b_1a_1-\bar b_{-1}a_{-1}$; thus $\delta\ne0$. If $a_0,b_0\ne0$, then $t_0$ is invertible and (2) gives $d=0$. If $a_0=0\ne b_0$, the two coordinates of $u_{12}=0$ give $a_{-1}=0$ and $a_1=0$, so $\delta=0$. If $b_0=0\ne a_0$, they give $b_1=0$ and $b_{-1}=0$, so $\delta=0$. Each case is contradictory unless $a_0=b_0=0$, i.e.\ $t_0=0$. The converse is part (1).

		(4) By Proposition~\ref{prop:red}(i), $T$ is binormal $\iff$ $A$ and $B$ are. For a complex $2\times2$ Toeplitz matrix the arguments of (1) and (2), with $\Hc$ replaced by $\C$, show: binormal $\iff$ normal or zero diagonal.
	\end{proof}
	
	\begin{corollary}\label{cor:cplx2}
		A complex $2\times 2$ Toeplitz matrix is normal $\iff$ $|t_1|=|t_{-1}|$, and it is binormal $\iff$ $|t_1|=|t_{-1}|$ or $t_0=0$.
	\end{corollary}
	
	\begin{proof}
		For complex entries and the usual conjugation, $d=|t_1|^2-|t_{-1}|^2$ and every non-zero $t_0$ is invertible. Apply Theorem~\ref{thm:n2}(1),(2).
	\end{proof}

	\begin{example}\label{ex:zd}
		Let
		\[
		T=\begin{pmatrix}\tfrac12(1+\mathbf{i}-\mathbf{j}-\mathbf{k})&1+\tfrac12 \mathbf{i}+\tfrac12 \mathbf{k}\\[2pt]\tfrac12(\mathbf{i}-\mathbf{k})&\tfrac12(1+\mathbf{i}-\mathbf{j}-\mathbf{k})\end{pmatrix}.
		\]
		In coordinates (using $q=\frac{a+b}2+\mathbf{j}\frac{a-b}2$), $t_0=(0,1+\mathbf{i})$, $t_1=(0,\mathbf{i})$, $t_{-1}=(1+\mathbf{i},1)$, so $t_0$ is a non-zero zero divisor.
		
		\emph{$\kappa=1$.} The components are $A=\begin{pmatrix}0&1+\mathbf{i}\\0&0\end{pmatrix}$ and $B=\begin{pmatrix}1+\mathbf{i}&1\\ \mathbf{i}&1+\mathbf{i}\end{pmatrix}$. $A$ has zero diagonal, hence is binormal, and is not normal since $|a_1|=0\ne\sqrt2=|a_{-1}|$; $B$ is normal by Corollary~\ref{cor:cplx2} since $|b_1|=|b_{-1}|=1$. By Theorem~\ref{thm:n2}(4), $T$ is binormal and not normal.
		
		\emph{$\kappa=2,3$.} By \eqref{eq:N}, $\Nm(t_1)=0$ (since $a_1=0$) while the first coordinate of $\Nm(t_{-1})$ equals $a_{-1}b_{-1}=1+\mathbf{i}$, resp.\ $\bar b_{-1}a_{-1}=1+\mathbf{i}$. Hence $d\ne0$ and $T$ is not normal; as $t_0\ne0$, Theorem~\ref{thm:n2}(3) shows that $T$ is not binormal.
		
		So the same matrix is binormal and not normal for the first conjugation and not even binormal for the other two. No complex $2\times 2$ Toeplitz matrix with $t_0\ne0$ is binormal and not normal (Corollary~\ref{cor:cplx2}).
	\end{example}
	
	\begin{example}\label{ex:neg}
		Let $T=\begin{pmatrix}0&\mathbf{j}\\1+\mathbf{k}&0\end{pmatrix}$, so $t_0=0$, $t_1=1+\mathbf{k}$, $t_{-1}=\mathbf{j}$. By Theorem~\ref{thm:n2}(1), $T$ is binormal for all three conjugations. For $\kappa=2$: $\Nm(t_1)=(1-\mathbf{k})(1+\mathbf{k})=2$ and $\Nm(t_{-1})=-\mathbf{j}^2=-1$, so
		\[
		T^\dagger T=\mathrm{diag}(2,-1),\qquad TT^\dagger=\mathrm{diag}(-1,2),
		\]
		which commute and differ. For $\kappa=3$: $\Nm(t_1)=(1+\mathbf{k})^2=2\mathbf{k}$ and $\Nm(t_{-1})=-1$, so $T^\dagger T=\mathrm{diag}(2\mathbf{k},-1)\ne\mathrm{diag}(-1,2\mathbf{k})=TT^\dagger$. For $\kappa=1$ the coordinates of $1+\mathbf{k}$ and $\mathbf{j}$ are $(1+\mathbf{i},1-\mathbf{i})$ and $(1,-1)$, so $\Nm(t_1)=2\ne1=\Nm(t_{-1})$. Thus $T$ is binormal and not normal for each conjugation, and for $\kappa=2$ the matrix $T^\dagger T$ has a negative diagonal entry, which is impossible for complex matrices with the usual adjoint. (It does occur if the transpose replaces the adjoint, so this phenomenon reflects the indefinite form $\Nm$ and not the zero divisors.)
	\end{example}
	
	\begin{example}\label{ex:Phi}
		Let $T=\begin{pmatrix}\mathbf{i}&\tfrac54\mathbf{i}+\tfrac34\mathbf{k}\\ \mathbf{i}&\mathbf{i}\end{pmatrix}$. Since $\frac54\mathbf{i}+\frac34\mathbf{k}=\mathbf{i}(\frac54+\frac34\mathbf{j})$, the coordinates are $t_0=t_1=(\mathbf{i},\mathbf{i})$ and $t_{-1}=(2\mathbf{i},\frac{\mathbf{i}}2)$, so $A=\mathbf{i}\begin{pmatrix}1&2\\1&1\end{pmatrix}$ and $B=\mathbf{i}\begin{pmatrix}1&\frac12\\1&1\end{pmatrix}$. For $\kappa=2$, $d=a_1b_1-a_{-1}b_{-1}=-1+1=0$, and for $\kappa=3$ the first coordinate of $d$ is $\bar b_1a_1-\bar b_{-1}a_{-1}=1-1=0$; so $T$ is normal for $\kappa=2,3$. For $\kappa=1$, $t_0=\mathbf{i}$ is invertible and $d=(1-4,\,1-\tfrac14)\ne0$, so $T$ is not binormal by Theorem~\ref{thm:n2}(2). Since $\Phi(T)$ is unitarily similar to $A\oplus B$ and $A$ is not binormal, $\Phi(T)$ is not binormal although $T$ is normal for $\kappa=2$: binormality of the complex adjoint is not equivalent to binormality of $T$ for $\kappa=2,3$.
	\end{example}
	
	\begin{example}\label{ex:isotropic}
		Let $T=\begin{pmatrix}0&0\\ \tfrac12(1+\mathbf{i}+\mathbf{j}+\mathbf{k})&0\end{pmatrix}$, whose entry $t_1$ has coordinates $(1+\mathbf{i},0)$. For $\kappa=2,3$, $\Nm(t_1)=0$ by \eqref{eq:N}, so $T^\dagger T=TT^\dagger=0$ although $T\ne0$, and $T$ is normal. For $\kappa=1$, $\Nm(t_1)=(2,0)\ne0$, so $d\ne0$ and $T$ is binormal (as $t_0=0$) but not normal. Over the complex numbers with the usual adjoint $T^*T=0$ forces $T=0$. With the transpose, $T^tT=0$ with $T\ne0$ already occurs for complex matrices, so this phenomenon reflects the indefinite form $\Nm$; what is special to $\Hc$ is that it happens for a matrix with a single non-zero entry, which is a zero divisor.
	\end{example}
	
	\subsection{Bidiagonal Matrices and Weighted Cyclic Shifts}
	
	Let $L\in M_n(\R)$ be the lower shift, $L_{i+1,i}=1$, and let $E_{ij}$ denote the matrix units. In this subsection $n\ge2$ is arbitrary.
	
	\begin{theorem}\label{thm:bidiag}
		Let $T=t_0I+t_1L\in M_n(\Hc)$ and put $\delta=\Nm(t_1)\in\Fix_\kappa$.
		\begin{enumerate}
			\item[(1)] $\Delta=\delta\,(E_{11}-E_{nn})$. $T$ is normal $\iff$ $\delta=0$, and binormal $\iff$ $\delta\cdot\bq t_1t_0=0$.
			\item[(2)] If $\kappa\in\{2,3\}$: $T$ is binormal and not normal $\iff$ $t_0=0$ and $\delta\ne0$. Moreover $T$ is normal $\iff$ $\delta=0$ $\iff$ $t_1$ is a zero divisor or $0$.
			\item[(3)] If $\kappa=1$: $T$ is normal $\iff$ $t_1=0$, and $T$ is binormal $\iff$ for each sign $\varepsilon$ one has $e_\varepsilon t_1=0$ or $e_\varepsilon t_0=0$.
		\end{enumerate}
	\end{theorem}
	
	\begin{proof}
		(1) Since $L$ is real, $T^\dagger=\bq t_0I+\bq t_1L^t$, and
		\[
		U=\Nm(t_0)I+\bq t_0t_1L+\bq t_1t_0L^t+\Nm(t_1)L^tL,\qquad V=\Nm(t_0)I+\bq t_0t_1L+\bq t_1t_0L^t+\Nm(t_1)LL^t .
		\]
		As $L^tL=\mathrm{diag}(1,\dots,1,0)$ and $LL^t=\mathrm{diag}(0,1,\dots,1)$, $\Delta=\delta(E_{11}-E_{nn})$, and $T$ is normal $\iff$ $\delta=0$. The matrix $U$ is tridiagonal with $u_{i,i+1}=\bq t_1t_0$. By Lemma~\ref{lem:crit}, $T$ is binormal $\iff$ $M=\Delta U$ is selfadjoint. Row $1$ of $M$ is $\delta$ times row $1$ of $U$, row $n$ is $-\delta$ times row $n$ of $U$, and the other rows vanish. Hence selfadjointness requires $\delta u_{1j}=\delta u_{nj}=0$ for $1<j<n$ and $2\delta u_{1n}=0$, the diagonal entries imposing nothing. The only possibly non-zero entries involved are $u_{12}=\bq t_1t_0$ and $u_{n,n-1}=\bq t_0t_1$ (and $u_{1n}=u_{12}$ if $n=2$); since $\bq\delta=\delta$ the second condition is the conjugate of the first, and everything reduces to $\delta\,\bq t_1t_0=0$.

		(2) If $\delta\ne0$, it is invertible in the field $\Fix_\kappa$, so $\bq t_1t_0=0$; also $\delta=\bq t_1t_1$ is invertible, so $t_1$ is invertible and $t_0=0$. Conversely, $t_0=0$ gives binormality by (1), and $T$ is not normal if $\delta\ne0$. By \eqref{eq:N}, $\delta=0\iff a_1b_1=0\iff t_1$ is a zero divisor or zero.

		(3) By \eqref{eq:N}, $\delta=(|a_1|^2,|b_1|^2)$ and $\bq t_1t_0=(\bar a_1a_0,\bar b_1b_0)$, so $\delta\bq t_1t_0=0$ $\iff$ ($a_1=0$ or $a_0=0$) and ($b_1=0$ or $b_0=0$). Finally $\delta=0\iff t_1=0$.
	\end{proof}

	In particular, a complex lower bidiagonal Toeplitz matrix $t_0I+t_1L$ is binormal if and only if $t_1=0$ or $t_0=0$, and normal if and only if $t_1=0$ (part (3) with both components equal).
	
	\begin{remark}\label{rem:monomial}
		The next statement is not specific to Toeplitz matrices: if $M\in M_n(\Hc)$ has at most one non-zero entry in each row and column (a monomial matrix), then $M^\dagger M$ and $MM^\dagger$ are diagonal, hence $M$ is binormal for every $\kappa$. A Toeplitz matrix is monomial if and only if either one diagonal is non-zero, or exactly the two diagonals $t_m$, $t_{m-n}$ ($0<m<n$) are non-zero, since two non-zero diagonals $t_m,t_l$ share a row as soon as $|m-l|<n$.
	\end{remark}

	\begin{proposition}\label{thm:qcirc}
		For $n\ge2$ and $s,q\in\Hc$ let $T_{s,q}=s\,L+q\,E_{1n}$, i.e.\ the Toeplitz matrix with $t_1=s$, $t_{-(n-1)}=q$ and all other $t_m=0$. Then, for each $\kappa$,
		\[
		T_{s,q}^\dagger T_{s,q}=\mathrm{diag}\big(\Nm(s),\dots,\Nm(s),\Nm(q)\big),\qquad T_{s,q}T_{s,q}^\dagger=\mathrm{diag}\big(\Nm(q),\Nm(s),\dots,\Nm(s)\big).
		\]
		Consequently $T_{s,q}$ is binormal, and it is normal $\iff$ $\Nm(q)=\Nm(s)$.
	\end{proposition}
	
	\begin{proof}
		The columns of $T_{s,q}$ are $s\,e_{j+1}$ for $1\le j<n$ and $q\,e_1$ for $j=n$, with pairwise disjoint supports, so $T^\dagger T$ is diagonal with the stated entries; likewise the rows are $q\,e_n^t$ and $s\,e_{i-1}^t$ ($i\ge2$), again with disjoint supports, so $TT^\dagger$ is diagonal with the stated entries. Two diagonal matrices commute, and they are equal if and only if the entries in position $1$ (and in position $n\ne1$) agree, i.e.\ $\Nm(q)=\Nm(s)$.
	\end{proof}
	
	\begin{example}\label{ex:cyc}
		(a) Let $s=\mathbf{i}$, $q=\mathbf{j}$ and $n=3$. By \eqref{eq:Nijk}, $\Nm(\mathbf{i})=(1,-1,1)$ and $\Nm(\mathbf{j})=(1,-1,-1)$ for $\kappa=1,2,3$. By Proposition~\ref{thm:qcirc}, $T$ is binormal for all three conjugations, normal for $\kappa=1,2$ (where $\Nm(\mathbf{i})=\Nm(\mathbf{j})$), and not normal for $\kappa=3$: $T^\dagger T=\mathrm{diag}(1,1,-1)\ne\mathrm{diag}(-1,1,1)=TT^\dagger$.

		(b) Let $s=\mathbf{k}$ and $q=1+\mathbf{i}+\mathbf{j}$. Then $\bq q=1-\mathbf{i}+\mathbf{j},\ 1+\mathbf{i}-\mathbf{j},\ 1-\mathbf{i}-\mathbf{j}$ for $\kappa=1,2,3$, so $\Nm(q)=3+2\mathbf{j},\ -1+2\mathbf{i},\ 1-2\mathbf{k}$, whereas $\Nm(\mathbf{k})=1,1,-1$. Hence $T$ is binormal and not normal for all three conjugations, and the last diagonal entry of $T^\dagger T$ lies in $\Fix_1,\Fix_2,\Fix_3$ respectively.
	\end{example}
	
	\begin{example}\label{ex:sum}
		Let $T_1=\mathbf{j}I$ and $T_2=\mathbf{i}L$ be $3\times3$ matrices. $T_1$ is normal, and $T_2$ is bidiagonal with $t_0=0$, hence binormal by Theorem~\ref{thm:bidiag}. Their sum $T=\mathbf{j}I+\mathbf{i}L$ has $t_0=\mathbf{j}$, $t_1=\mathbf{i}$, so $\delta=\Nm(\mathbf{i})=\pm1$ is invertible, and $\bq t_1t_0=\bar{\mathbf{i}}\,\mathbf{j}=-\mathbf{k},\mathbf{k},-\mathbf{k}$ for $\kappa=1,2,3$ is non-zero. By Theorem~\ref{thm:bidiag}(1), $T$ is not binormal for any $\kappa$.
	\end{example}
	
	\subsection{The $3\times 3$ Case}
	
	For $T=(t_{i-j})_{i,j=1}^3$, Lemma~\ref{lem:uv} gives
	\begin{align*}
		u_{11}&=\Nm(t_0)+\Nm(t_1)+\Nm(t_2), & u_{12}&=\bq t_0t_{-1}+\bq t_1t_0+\bq t_2t_1,\\
		u_{22}&=\Nm(t_{-1})+\Nm(t_0)+\Nm(t_1), & u_{13}&=\bq t_0t_{-2}+\bq t_1t_{-1}+\bq t_2t_0,\\
		u_{33}&=\Nm(t_{-2})+\Nm(t_{-1})+\Nm(t_0), & u_{23}&=\bq t_{-1}t_{-2}+\bq t_0t_{-1}+\bq t_1t_0,
	\end{align*}
	with $u_{ji}=\bq u_{ij}$. Put
	\[
	\delta:=u_{11}-u_{33}\in\Fix_\kappa,\qquad \varepsilon:=u_{12}-u_{23}.
	\]
	By Lemma~\ref{lem:V}, $V=\begin{pmatrix}u_{33}&u_{23}&u_{13}\\ \bq u_{23}&u_{22}&u_{12}\\ \bq u_{13}&\bq u_{12}&u_{11}\end{pmatrix}$, so
	\[
	\Delta=U-V=\begin{pmatrix}\delta&\varepsilon&0\\ \bq\varepsilon&0&-\varepsilon\\ 0&-\bq\varepsilon&-\delta\end{pmatrix}
	\]
	(the entry $(1,3)$ is $u_{13}-u_{13}=0$ and the entry $(2,2)$ is $u_{22}-u_{22}=0$). In particular $T$ is normal $\iff$ $\delta=0$ and $\varepsilon=0$.

	We now give a complete description in terms of the entries of $T$. By Proposition~\ref{prop:pair} it suffices to treat a pair $(A,C)$ of complex $3\times 3$ Toeplitz matrices, $A=(a_{i-j})$, $C=(c_{i-j})$, and the condition $[CA,AC]=0$. Put
	\[
	\begin{aligned}
		\delta&=a_1c_{-1}+a_2c_{-2}-a_{-1}c_1-a_{-2}c_2,\\
		\varepsilon&=a_1c_{-2}-a_{-2}c_1,\qquad \varepsilon'=a_2c_{-1}-a_{-1}c_2,\\
		P&=2a_0c_0+a_1c_{-1}+a_2c_{-2}+a_{-1}c_1+a_{-2}c_2,\\
		S&=2a_0c_0+2a_1c_{-1}+2a_{-1}c_1,\\
		Q&=2a_0c_{-1}+2a_{-1}c_0+a_1c_{-2}+a_{-2}c_1,\\
		Q'&=2a_0c_1+2a_1c_0+a_2c_{-1}+a_{-1}c_2,\\
		R&=2(a_0c_{-2}+a_{-1}c_{-1}+a_{-2}c_0),\\
		R'&=2(a_0c_2+a_1c_1+a_2c_0),
	\end{aligned}
	\]
	and the vectors
	\[
	v=\big(\delta\varepsilon,\ \delta\varepsilon',\ -\varepsilon^2,\ -\varepsilon'^2,\ \delta^2-\varepsilon\varepsilon'\big),\qquad w=\big(Q,\ Q',\ R,\ R',\ P-S\big)\in\C^5 .
	\]
	
	\begin{theorem}\label{thm:n3full}
		Let $A,C$ be complex $3\times3$ Toeplitz matrices and $D=CA-AC$, $\Sigma_0=CA+AC$.
		\begin{enumerate}
			\item[(a)] $D=\begin{pmatrix}\delta&\varepsilon&0\\ \varepsilon'&0&-\varepsilon\\ 0&-\varepsilon'&-\delta\end{pmatrix}$, and $CA=AC$ if and only if $v=0$.
			\item[(b)] $[CA,AC]=0$ if and only if $w\wedge v=0$, i.e.\ if and only if $v=0$ or $w=\gamma v$ for some $\gamma\in\C$. Equivalently: $D=0$, or $CA+AC=\alpha I+\gamma D^2$ for some $\alpha,\gamma\in\C$.
		\end{enumerate}
		If $C=A^*$ (the case $\kappa=1$), then $\gamma$ is real. Consequently a $3\times3$ Toeplitz matrix $T=e_+A+e_-B$ over $\Hc$ is binormal if and only if the condition $w\wedge v=0$ holds for the pair(s) of Proposition~\ref{prop:pair}, i.e.\ for $(A,A^*)$ and $(B,B^*)$ if $\kappa=1$, for $(A,B^t)$ if $\kappa=2$, and for $(A,B^*)$ if $\kappa=3$; it is normal if and only if $v=0$ for these pairs.
	\end{theorem}
	
	\begin{proof}
		Put $X=CA$, $Y=AC$. For a Toeplitz matrix $G$ one has $JG^tJ=G$ (Section 3), so $JX^tJ=JA^tC^tJ=(JA^tJ)(JC^tJ)=AC=Y$. Hence $\Sigma_0=X+Y$ is persymmetric ($J\Sigma_0^tJ=\Sigma_0$) and $D=X-Y$ is skew-persymmetric ($JD^tJ=-D$), which gives the displayed shape of $D$ and the shape $\Sigma_0=\begin{pmatrix}P&Q&R\\ Q'&S&Q\\ R'&Q'&P\end{pmatrix}$; direct multiplication gives the displayed formulas. Also $[X,Y]=-\tfrac12[\Sigma_0,D]$. Part (a): $D=0\iff\delta=\varepsilon=\varepsilon'=0\iff v=0$ (if $v=0$ then $\varepsilon^2=\varepsilon'^2=0$ and then $\delta^2=0$).
		
		(b) If $D=0$ both sides hold. Let $D\ne0$. Its characteristic polynomial is $x^3-\rho^2x$ with $\rho^2=\delta^2+2\varepsilon\varepsilon'$, and $D^2=\begin{pmatrix}\delta^2+\varepsilon\varepsilon'&\delta\varepsilon&-\varepsilon^2\\ \delta\varepsilon'&2\varepsilon\varepsilon'&\delta\varepsilon\\ -\varepsilon'^2&\delta\varepsilon'&\delta^2+\varepsilon\varepsilon'\end{pmatrix}$. We claim that $D$ is non-derogatory, i.e.\ $I,D,D^2$ are linearly independent. If $\rho\ne0$, $D$ has three distinct eigenvalues $0,\pm\rho$. If $\rho=0$, $D$ is nilpotent and non-zero; a rank-one matrix would have all $2\times2$ minors zero, and the minors $-\varepsilon\varepsilon'$, $-\varepsilon^2$, $-\varepsilon'^2$ of $D$ would force $\varepsilon=\varepsilon'=0$ and then $\delta=0$, i.e.\ $D=0$; so $D$ has rank two, hence a single Jordan block and $D^2\ne0$. A non-derogatory matrix has centralizer $\mathrm{span}\{I,D,D^2\}$. Therefore $[\Sigma_0,D]=0$ if and only if $\Sigma_0=\alpha I+\beta D+\gamma D^2$. The map $\pi(G)=JG^tJ$ is linear, fixes $I$, reverses products, and $\pi(D)=-D$, $\pi(D^2)=D^2$; applying it and using $\pi(\Sigma_0)=\Sigma_0$ and the independence of $I,D,D^2$ gives $\beta=0$. Comparing entries of $\Sigma_0=\alpha I+\gamma D^2$ gives $Q=\gamma\delta\varepsilon$, $Q'=\gamma\delta\varepsilon'$, $R=-\gamma\varepsilon^2$, $R'=-\gamma\varepsilon'^2$, $P-S=\gamma(\delta^2-\varepsilon\varepsilon')$, i.e.\ $w=\gamma v$; conversely these five equations make $\Sigma_0-\gamma D^2$ a matrix with zero off-diagonal part and equal diagonal entries, i.e.\ a scalar matrix. Since $v\ne0$ when $D\ne0$, $w\wedge v=0$ means $w=\gamma v$.
		
		For $C=A^*$ we have $c_m=\bar a_{-m}$, so $\delta,P,S$ are real, $\varepsilon'=\bar\varepsilon$, $Q'=\bar Q$, $R'=\bar R$. Complex conjugation exchanges the first and second, and the third and fourth, coordinates of both $w$ and $v$; so $w=\gamma v$ with $v\ne0$ gives $w=\bar\gamma v$, hence $\gamma=\bar\gamma$. The statement for $\Hc$ is Proposition~\ref{prop:pair}.
	\end{proof}
	
	Normal $3\times 3$ Toeplitz matrices are those with $v=0$; the binormal non-normal ones are exactly the solutions of $w=\gamma v$ with $v\ne0$. In contrast with the two types of normal Toeplitz matrices, this is a proportionality condition on quadratic expressions in the entries.
	
	\begin{example}\label{ex:n3new}
		Let $T=\begin{pmatrix}-1&-\mathbf{j}&2\\ -\mathbf{j}&-1&-\mathbf{j}\\ 2&-\mathbf{j}&-1\end{pmatrix}$, a symmetric Toeplitz matrix with $t_0=-1$, $t_{\pm1}=-\mathbf{j}$, $t_{\pm2}=2$. Its components are the real symmetric, hence normal, matrices $A=\begin{pmatrix}-1&-1&2\\ -1&-1&-1\\ 2&-1&-1\end{pmatrix}$ and $B=\begin{pmatrix}-1&1&2\\ 1&-1&1\\ 2&1&-1\end{pmatrix}$, so $T$ is normal for $\kappa=1$. For $\kappa=2$, and also $\kappa=3$ since $B$ is real, the pair is $(A,C)=(A,B)$ with $a_0=-1$, $a_{\pm1}=-1$, $a_{\pm2}=2$, $c_0=-1$, $c_{\pm1}=1$, $c_{\pm2}=2$. Then $\delta=0$, $\varepsilon=-4$, $\varepsilon'=4$, so $v=(0,0,-16,-16,16)$ and $w=(0,0,-10,-10,10)=\tfrac58v$. Hence $w\wedge v=0$ and $T$ is binormal by Theorem~\ref{thm:n3full}, but $v\ne0$, so $T$ is not normal.
	\end{example}
	
		\begin{corollary}\label{cor:diag}
		If $U$ is diagonal (the columns of $T$ are pairwise orthogonal), then $T$ is binormal.
	\end{corollary}
	
	\begin{proof}
		If $U$ is diagonal, then so is $V=J\bq UJ$ (Lemma~\ref{lem:V}); two diagonal matrices commute.
	\end{proof}

	\begin{example}\label{ex:mixed}
		Let $T=(t_{i-j})$ be the symmetric Toeplitz matrix with
		\[
		t_0=\tfrac12(1-\mathbf{i}+\mathbf{j}+\mathbf{k}),\qquad t_{\pm1}=1-\mathbf{i}-\mathbf{j}-\mathbf{k},\qquad t_{\pm2}=-1+\mathbf{i}-\mathbf{j}-\mathbf{k}.
		\]
		Its coordinates are $t_0=(1,-\mathbf{i})$, $t_{\pm1}=(-2\mathbf{i},2)$, $t_{\pm2}=(-2,2\mathbf{i})$, so the components are
		\[
		A=\begin{pmatrix}1&-2\mathbf{i}&-2\\-2\mathbf{i}&1&-2\mathbf{i}\\-2&-2\mathbf{i}&1\end{pmatrix},\qquad B=-\mathbf{i}\,\bar A .
		\]
		For $\kappa=1$ we examine $A$ and $B$ (Proposition~\ref{prop:red}(i)). Direct multiplication gives $A^*A=\begin{pmatrix}9&4\mathbf{i}&0\\-4\mathbf{i}&9&-4\mathbf{i}\\0&4\mathbf{i}&9\end{pmatrix}$, $AA^*=\overline{A^*A}$ and $A^*A\cdot AA^*=AA^*\cdot A^*A=\begin{pmatrix}65&0&-16\\0&49&0\\-16&0&65\end{pmatrix}$. So $A$ is binormal and not normal. Since $B=-\mathbf{i}\bar A$, $B^*B=\overline{A^*A}$ and $BB^*=\overline{AA^*}$, and complex conjugation preserves commutation, so $B$ is binormal. Hence $T$ is binormal and not normal for $\kappa=1$.
	\end{example}
	
	\begin{example}\label{ex:t0zero}
		Let $\kappa=1$ and $T=\begin{pmatrix}0&0&\mathbf{j}\\ \mathbf{i}&0&0\\ \mathbf{j}&\mathbf{i}&0\end{pmatrix}$, i.e.\ $t_0=0$, $t_1=\mathbf{i}$, $t_2=\mathbf{j}$, $t_{-1}=0$, $t_{-2}=\mathbf{j}$. The first component is $A=\begin{pmatrix}0&0&1\\ \mathbf{i}&0&0\\ 1&\mathbf{i}&0\end{pmatrix}$, and the pair is $(A,A^*)$, $c_m=\bar a_{-m}$. Theorem~\ref{thm:n3full} gives $\delta=1$, $\varepsilon=\mathbf{i}$, $\varepsilon'=-\mathbf{i}$, so $v=(\mathbf{i},-\mathbf{i},1,1,0)$, while $Q=\mathbf{i}$, $Q'=-\mathbf{i}$, $R=R'=0$, $P-S=1$ give $w=(\mathbf{i},-\mathbf{i},0,0,1)$. These vectors are not parallel, so $T$ is not binormal although $t_0=0$. In contrast, for $2\times2$ matrices and for bidiagonal matrices of every size, $t_0=0$ forces binormality (Theorems~\ref{thm:n2} and~\ref{thm:bidiag}).
	\end{example}
	
	\subsection{The Shift Property}\label{sec:shift}
	
	In this subsection $\lambda\in\Hc$ is a scalar and $T\in M_n(\Hc)$ is an arbitrary square matrix; we use $U=T^\dagger T$, $V=TT^\dagger$, $\Delta=U-V$ as before.
	
	\begin{lemma}\label{lem:shift}
		For every $\lambda\in\Hc$, with the selfadjoint matrix $S:=\lambda T^\dagger+\bq\lambda T$,
		\[
		(T-\lambda I)^\dagger(T-\lambda I)=U-S+\Nm(\lambda)I,\qquad (T-\lambda I)(T-\lambda I)^\dagger=V-S+\Nm(\lambda)I .
		\]
		In particular the difference of the two products for $T-\lambda I$ equals $\Delta$, and $T$ is normal $\iff$ $T-\lambda I$ is normal.
	\end{lemma}
	
	\begin{proof}
		Since $(T-\lambda I)^\dagger=T^\dagger-\bq\lambda I$ and scalars are central, $(T^\dagger-\bq\lambda)(T-\lambda)=U-\lambda T^\dagger-\bq\lambda T+\Nm(\lambda)I$ and $(T-\lambda)(T^\dagger-\bq\lambda)=V-\bq\lambda T-\lambda T^\dagger+\Nm(\lambda)I$. Moreover $S^\dagger=S$, and subtracting gives $\Delta$.
	\end{proof}
	
	\begin{lemma}\label{lem:nilp}
		Let $K\in M_n(\Hc)$ satisfy $\tr(K^m)=0$ for all $m\ge1$, where $\tr$ denotes the trace of matrix. Then $K$ is nilpotent.
	\end{lemma}
	
	\begin{proof}
		Writing $K=e_+K_A+e_-K_B$. By the ring isomorphism of Section~2, $\tr(K^m)=e_+\tr(K_A^m)+e_-\tr(K_B^m)$, so $\tr(K_A^m)=\tr(K_B^m)=0$ for all $m\ge1$. Over $\C$ the vanishing of the power sums $p_1,\dots,p_n$ of the eigenvalues implies, by Newton's identities in characteristic zero, that all elementary symmetric functions vanish; hence the characteristic polynomials are $x^n$ and $K_A^n=K_B^n=0$ by Cayley--Hamilton, so $K^n=0$.
	\end{proof}
	
	\begin{theorem}\label{thm:nilp}
		Let $T\in M_n(\Hc)$ and $\lambda\in\Hc$, and assume that $T$ and $T-\lambda I$ are both binormal. If there is $\mu\in\Hc$ for which
		\[
		c:=\lambda\bq\mu-\bq\lambda\mu
		\]
		is invertible, then $\Delta=T^\dagger T-TT^\dagger$ is nilpotent. Such a $\mu$ exists
		\begin{itemize}
			\item for $\kappa=2$ and $\kappa=3$ and every $\lambda\ne0$;
			\item for $\kappa=1$ if and only if $\lambda$ is not a zero divisor.
		\end{itemize}
	\end{theorem}
	
	\begin{proof}
		\emph{Step 1.} Binormality of $T$ means $[U,\Delta]=0$ (as $V=U-\Delta$); by Lemma~\ref{lem:shift}, binormality of $T-\lambda I$ means $[U-S,\Delta]=0$. Subtracting, $[S,\Delta]=0$.

		\emph{Step 2.} Put $W=\mu T^\dagger+\bq\mu T$. Expanding bilinearly (scalars are central) and using $[T^\dagger,T^\dagger]=[T,T]=0$, $[T^\dagger,T]=\Delta=-[T,T^\dagger]$,
		\[
		[S,W]=(\lambda\bq\mu-\bq\lambda\mu)\,\Delta=c\,\Delta,\qquad\text{so}\quad \Delta=c^{-1}[S,W].
		\]

		\emph{Step 3.} For $m\ge1$, using $\Delta^{m-1}S=S\Delta^{m-1}$,
		\[
		\Delta^m=c^{-1}\Delta^{m-1}(SW-WS)=c^{-1}[S,\Delta^{m-1}W].
		\]
		Over a commutative ring $\tr(XY)=\tr(YX)$, so $\tr(\Delta^m)=0$ for all $m\ge1$, and $\Delta$ is nilpotent by Lemma~\ref{lem:nilp}.

		\emph{Step 4:} We write $\lambda=(\lambda_1,\lambda_2)$, $\mu=(\mu_1,\mu_2)$. For $\kappa=1$, $c=(\lambda_1\bar\mu_1-\bar\lambda_1\mu_1,\ \lambda_2\bar\mu_2-\bar\lambda_2\mu_2)$. If $\lambda_1,\lambda_2\ne0$, take $\mu_r=\mathrm{i}\lambda_r$; then each coordinate of $c$ equals $-2\mathrm{i}|\lambda_r|^2\ne0$. If $\lambda_1=0$ (or $\lambda_2=0$), the corresponding coordinate of $c$ vanishes for every $\mu$, so $c$ is never invertible. For $\kappa=2$, $c=(\lambda_1\mu_2-\lambda_2\mu_1,\ \lambda_2\mu_1-\lambda_1\mu_2)$; if $\lambda_1\ne0$ take $\mu=(0,1)$, giving $c=(\lambda_1,-\lambda_1)$, and if $\lambda_1=0\ne\lambda_2$ take $\mu=(1,0)$, giving $c=(-\lambda_2,\lambda_2)$. For $\kappa=3$ apply the case $\kappa=2$ to $\tau(\lambda)\ne0$ (Lemma~\ref{lem:tau}): $\tau$ carries the element $c$ for $\kappa=2$ to the element $c$ for $\kappa=3$ and preserves invertibility.
	\end{proof}
	
	\begin{lemma}\label{lem:cshift}
		Let $A\in M_n(\C)$ and $0\ne\lambda\in\C$. If $A$ and $A-\lambda I$ are both binormal, then $A$ is normal.
	\end{lemma}
	
	\begin{proof}
		The proof of Theorem~\ref{thm:nilp} applies exactly to complex matrices with the usual adjoint: with $\mu=\mathrm{i}\lambda$ one gets $c=\lambda\bar\mu-\bar\lambda\mu=-2\mathrm{i}|\lambda|^2\ne0$, so $\Delta=A^*A-AA^*$ is nilpotent. But $\Delta$ is Hermitian, hence its eigenvalues are real, and nilpotency forces them all to vanish; so $\Delta=0$.
	\end{proof}
	This is the classical shift property of Campbell \cite{C2}; the trace argument used here goes back to Jacobson \cite{Ja}.
	
	\begin{theorem}\label{thm:shift1}
		Let $\kappa=1$, let $T=e_+A+e_-B$ and $T-\lambda I$ be binormal, where $\lambda=(\lambda_+,\lambda_-)$ in coordinates. Then $A$ is normal if $\lambda_+\ne0$, and $B$ is normal if $\lambda_-\ne0$. In particular $T$ is normal whenever $\lambda$ is not a zero divisor, for instance for every non-zero real or complex $\lambda$.
	\end{theorem}
	
	\begin{proof}
		Since $T-\lambda I=e_+(A-\lambda_+I)+e_-(B-\lambda_-I)$, Proposition~\ref{prop:red}(i) shows that $A$, $A-\lambda_+I$, $B$, $B-\lambda_-I$ are binormal. Apply Lemma~\ref{lem:cshift} to $A$ if $\lambda_+\ne0$ and to $B$ if $\lambda_-\ne0$; if $\lambda$ is not a zero divisor both components are normal, hence so is $T$.
	\end{proof}
	
	\begin{corollary}\label{cor:diagDelta}
		Let $\kappa\in\{2,3\}$ and let $0\ne\lambda\in\Hc$ be arbitrary, zero divisors included. Suppose $\Delta=T^\dagger T-TT^\dagger$ is a diagonal matrix. If $T$ and $T-\lambda I$ are binormal, then $T$ is normal. This applies to all $2\times 2$ Toeplitz matrices, all lower bidiagonal Toeplitz matrices (Theorem~\ref{thm:bidiag}(1)) and all weighted cyclic shifts (Proposition~\ref{thm:qcirc}).
	\end{corollary}
	
	\begin{proof}
		By Theorem~\ref{thm:nilp}, $\Delta$ is nilpotent. A diagonal matrix $\mathrm{diag}(x_1,\dots,x_n)$ with $x_r=(a_r,b_r)$ is nilpotent $\iff$ $x_r^m=0$ for some $m$, i.e.\ $a_r^m=b_r^m=0$, i.e.\ $x_r=0$. So $\Delta=0$. For the listed classes, $\Delta$ is diagonal by Theorem~\ref{thm:n2}(1) (where $\Delta=\mathrm{diag}(d,-d)$), Theorem~\ref{thm:bidiag}(1), and Proposition~\ref{thm:qcirc} (where $\Delta=\mathrm{diag}(\Nm(s)-\Nm(q),0,\dots,0,\Nm(q)-\Nm(s))$).
	\end{proof}
	
	Thus, within these classes, for the second and third conjugations the shift property holds for every non-zero scalar, zero divisors included, whereas for the first conjugation it fails for zero-divisor scalars. This is only a statement about matrices with diagonal $\Delta$; for general matrices the situation is the opposite (Example~\ref{ex:genfail}).
	
	\begin{example}\label{ex:shiftfail}
		Let $\kappa=1$, $T=\tfrac12(1+\mathbf{i}+\mathbf{j}+\mathbf{k})\,L$ (of size $3$) and $\lambda=\mathbf{i}-\mathbf{k}=\mathbf{i}(1-\mathbf{j})$. In coordinates $t_1=(1+\mathbf{i},0)$ and $\lambda=(0,2\mathbf{i})$. Then $t_0=0$, so $T$ is binormal by Theorem~\ref{thm:bidiag}, and it is not normal since $t_1\ne0$. The matrix $T-\lambda I$ has $t_0=-\lambda=(0,-2\mathbf{i})$ and the same $t_1$; by Theorem~\ref{thm:bidiag}(3) it is binormal, since $a_0=0$ and $b_1=0$. So $T$ and $T-\lambda I$ are binormal, $\lambda\ne0$, and $T$ is not normal. This does not contradict Theorem~\ref{thm:shift1}, because $\lambda$ is a zero divisor. For $\kappa=2,3$ the same matrix is normal, in accordance with Corollary~\ref{cor:diagDelta}, because $\Nm(t_1)=0$.
	\end{example}
	The following example shows that the shift property fails for every $\lambda$ when $\kappa=2,3$.
	\begin{example}\label{ex:genfail}
		 Let $T=e_+E_{23}+e_-E_{21}\in M_3(\Hc)$, $e_\pm=\frac12(1\pm\mathbf{j})$. Its components are $A=E_{23}$, $B=E_{21}$, and for $\kappa=2,3$ the matrix $C$ of Proposition~\ref{prop:pair} is $C=E_{12}$. Then $CA=E_{13}$ and $AC=0$, so $T$ is not normal, while $[CA,AC]=0$, so $T$ is binormal. Let $\lambda=(\lambda_1,\lambda_2)\in\Hc$ be arbitrary and put $\mu=\lambda_2$ ($\kappa=2$), $\mu=\bar\lambda_2$ ($\kappa=3$). For $T-\lambda I$ the pair is $(A-\lambda_1I,\,C-\mu I)$ and
		\[
		(C-\mu)(A-\lambda_1)=E_{13}-\lambda_1E_{12}-\mu E_{23}+\lambda_1\mu I,\qquad (A-\lambda_1)(C-\mu)=-\mu E_{23}-\lambda_1E_{12}+\lambda_1\mu I .
		\]
		They differ by $E_{13}$, which commutes with $E_{12}$ and $E_{23}$, so the two products commute: $T-\lambda I$ is binormal for \emph{every} $\lambda\in\Hc$ (zero divisors included), and never normal. Here $\Delta$ is nilpotent and non-zero, so Theorem~\ref{thm:nilp} cannot be improved. For $\kappa=1$ the component $A=E_{23}$ is binormal, but $A-\lambda_1I$ is not for $\lambda_1\ne0$ (Lemma~\ref{lem:cshift}), in accordance with Theorem~\ref{thm:shift1}. The matrix is not Toeplitz.
	\end{example}

\medskip
	 Let $\kappa\in\{2,3\}$, let $T$ be a Toeplitz matrix of size at least three, and suppose that $T$ and $T-\lambda I$ are binormal for some $\lambda\ne0$. Is $T$ normal? By Theorem~\ref{thm:nilp} the matrix $\Delta$ is nilpotent. For matrices that are not Toeplitz the answer is negative (Example~\ref{ex:genfail}), and for $2\times2$ Toeplitz matrices it is affirmative (Corollary~\ref{cor:diagDelta}). For $3\times3$ matrices Theorem~\ref{thm:n3full} reduces the question to the pair formulation: one needs $[CA,AC]=0$ and $[D,\lambda_1C+\mu A]=0$ with $D=[C,A]$. At present we are not able to answer this question, and it will be addressed elsewhere.
	
		\section*{Acknowledgments}
	\section*{Declarations}
	\subsection*{Competing Interests}
	The author declare that he has no known competing financial interests or personal relationships that could have appeared to influence the work reported in this article.
	\subsection*{Funding Information}
	This research did not receive any specific grant or funding.
	\subsection*{Author Contribution}
	The author solely contributed to the preparation of the article. 
	\subsection*{Data Availability Statement}
	Data sharing is not applicable to this article as no data sets were generated or analyzed during the current study.

\end{document}